\documentclass[sn-mathphys-num,draft=false]{sn-jnl}
\usepackage{graphicx} 
\usepackage{amsmath}
\usepackage{fixmetodonotes}
\usepackage{amsfonts}
\usepackage{amssymb}
\usepackage{amsthm}
\usepackage{siunitx}
\usepackage{hyperref}
\usepackage{caption}
\usepackage{float}
\usepackage{subcaption}
\usepackage[normalem]{ulem}
\usepackage{cancel}
\usepackage{booktabs}
\usepackage{stmaryrd}
\usepackage{multirow}

\DeclareSIUnit{\Osm}{\text{Osm}}

\begin{document}

\title{Parameter estimation for evaporation-driven tear film model in two space dimensions}

\author[1]{{\fnm{Qinying} \sur{Chen}} \email{cqyzxy@udel.edu}}

\author[1]{{\fnm{Tobin A.} \sur{Driscoll}} \email{driscoll@udel.edu}}

\affil[1]{\orgdiv{Department of Mathematical Sciences}, \orgname{University of Delaware}, \orgaddress{\city{Newark}, \postcode{19716}, \state{DE}, \country{USA}}}

\abstract{The tear film (TF) plays a critical role in maintaining ocular surface health, and its disruption through tear breakup (TBU) is closely associated with dry eye disease. Evaporation-driven thinning is a primary mechanism underlying TBU, yet quantitative in vivo estimates of key physical parameters remain limited. In this work, we fit an evaporation-driven TF thinning model, originally developed by Braun et al. and extended to two dimensions using proper orthogonal decomposition (POD) by Chen et al., to experimental fluorescence (FL) imaging data from normal subjects. The use of dimension reduction enables efficient solution of the governing PDEs and facilitates parameter estimation from imaging data. Our results provide in vivo estimates of evaporation-related and thinning parameters within TBU regions. These findings enhance understanding of TF thinning and dry-spot formation and establish a quantitative baseline for comparison with dry eye patient data.}

\keywords{Tear film, Dry eye disease, Fluorescent imaging, Optimization}
\maketitle

\section{Introduction}
Each blink causes the upper eyelid to descend and rise, spreading a thin layer of fluid known as the tear film (TF) across the ocular surface \cite{Doane80}. In healthy eyes, the TF serves multiple essential functions: it lubricates the eye and eyelids, provides antimicrobial protection, maintains a smooth refractive surface, and delivers oxygen and nutrients to the avascular corneal epithelium \cite{lempDefinitionClassificationDry2007,willcoxTFOSDEWSII2017}. Disruption of this film, commonly referred to as tear breakup (TBU), exposes the ocular surface to harmful stimuli \cite{king2018mechanisms}. Persistent TBU and tear-film dysfunction are widely considered key contributors to the onset and progression of dry eye disease (DED). The prevalence of DED is substantial, though estimates vary depending on diagnostic definitions \cite{stapletonDEWSIIepi2017}. Beyond discomfort, DED significantly impairs visual quality and compromises ocular surface health \cite{nelsonTFOSDEWSII2017}. Among the major subtypes of dry eye, evaporative dry eye (EDE) is the most common \cite{lempDefinitionClassificationDry2007} and is primarily attributed to excessive tear loss due to evaporation \cite{linDryEyeDisease2014,OCEANreport2013}.

The tear film is a thin, multilayered liquid structure that rapidly reforms following each blink \cite{braunTearFilm2018}. It is commonly described as consisting of three distinct layers: an outer lipid layer, approximately 20--100 nm thick \cite{braunDynamicsFunctionTear2015}; a central aqueous layer, composed primarily of water and several microns in thickness \cite{hollyFormationRuptureTear1973}; and an inner mucin-rich layer, known as the glycocalyx, which is about half a micron thick and coats the ocular surface \cite{king-smithThicknessTearFilm2004}. The lipid layer plays a critical role in reducing evaporation from the tear film \cite{mishimaOilyLayerTear1961}, while an intact and healthy glycocalyx promotes smooth fluid transport along the ocular surface \cite{gipsonDistributionMucinsOcular2004}. Most of the aqueous component is supplied by the lacrimal gland, with secretion occurring predominantly near the temporal canthus \cite{darttNeuralRegulationLacrimal2009}. Additional water enters the tear film through osmotic transport from the ocular epithelia \cite{braunDynamicsTear2012}.

TBU is defined by the formation of localized dry spots on the ocular surface \cite{nornMICROPUNCTATEFLUORESCEINVITAL1970} and is frequently driven by evaporation \cite{lempDefinitionClassificationDry2007,willcoxTFOSDEWSII2017}. The tear breakup time (TBUT) measures the interval between a blink and the initial appearance of such a dry spot \cite{nornMICROPUNCTATEFLUORESCEINVITAL1970}. In clinical practice, TBUT is determined subjectively, relying on the clinician's visual assessment \cite{nornMICROPUNCTATEFLUORESCEINVITAL1970}, and may involve averaging estimates across multiple observers to improve reliability \cite{choReliabilityTearBreakup1992}. TBUT serves as a standard indicator of tear-film stability: shorter breakup times reflect diminished tear-film quality, whereas longer TBUT values indicate a more stable and robust tear film \cite{dibajniaTearFilmBreakup2012}.

Tear-breakup (TBU) models are typically formulated on short spatial domains to capture local tear-film (TF) dynamics while neglecting meniscus effects. Evaporation is the primary driving mechanism in many such models \cite{PengEtal2014,braunDynamicsFunctionTear2015,braunTearFilm2018}, with Peng et al. \cite{PengEtal2014} additionally incorporating osmolarity transport in the aqueous layer (AL) and osmosis across the AL--cornea interface, as well as spatially varying lipid layer (LL) thickness, to show that osmolarity diffusion prevents osmosis from arresting thinning as predicted by spatially uniform models \cite{braunDynamicsTear2012,braunDynamicsFunctionTear2015}. Simpler evaporation-driven models that include fluorescein transport were later used to interpret tear-film visualization experiments \cite{braunDynamicsFunctionTear2015,braunTearFilm2018}, and Zhong et al. \cite{zhongDynamicsFluorescentImaging2019} developed a one-dimensional PDE model combining evaporation and Marangoni effects. Incorporating fluorescein transport and fluorescence enabled model fitting to in vivo data within TBU regions to estimate otherwise inaccessible parameters \cite{lukeParameterEstimation2020,lukeParameterEstimationMixedMechanism2021}. More recently, spatially lumped ODE models have been fit to fluorescence data from small TBU spots and streaks \cite{lukeFittingSimplifiedModels2021}, and coupling these models with neural-network--based data extraction has greatly expanded the number of analyzable TBU instances \cite{driscollFittingODEModels2023}.

Imaging of the tear film is a crucial tool for analyzing its dynamics. Common imaging modalities include fluorescence (FL) imaging \cite{king-smithTearFilmInterferometry2014}, spectral interferometry \cite{king-smithApplicationNovelInterferometric2010}, and optical coherence tomography \cite{wangPrecornealPrePostlens2003}. The injection of dyes such as fluorescein has been used to stain epithelial cells \cite{nornMICROPUNCTATEFLUORESCEINVITAL1970}, estimate tear drainage rates or turnover times \cite{webberContinuousFluorophotometricMethod1986}, visualize overall tear-film dynamics \cite{benedettoVivoObservationTear1984,begleyQuantitativeAnalysisTear2013}, estimate tear-film breakup times \cite{nornMICROPUNCTATEFLUORESCEINVITAL1970}, and identify tear-breakup (TBU) regions. Simultaneous multimodal imaging has also been used to aid interpretation of tear-film dynamics \cite{himebaughScaleSpatialDistribution2012}.

Luke et al. \cite{lukeParameterEstimation2020} developed a parameter-estimation framework that fits fluorescence (FL) imaging data to evaporation-driven tear-film-thinning models \cite{braunTearFilm2018}, yielding physiologically realistic estimates of evaporation rates, dry-spot sizes, and thinning rates consistent with experimental measurements \cite{nicholsThinningRatePrecorneal2005a}. However, the computational cost of the PDE-based approach limited the analysis to a small number of thinning events. Chen et al. \cite{chenEvaporationdrivenTearFilm2024} introduced a dimension-reduction approach based on proper orthogonal decomposition (POD) to accelerate the solution of the governing PDEs, thereby significantly improving the computational efficiency of the associated inverse problem for parameter estimation. Driscoll et al. \cite{driscollFittingODEModels2023} introduced a simplified model with automated tear-breakup (TBU) detection, enabling the analysis of hundreds of thinning instances while preserving trends observed in PDE models, albeit at the expense of detailed spatial information. Together, these studies demonstrate that more efficient fitting of imaging data can substantially expand the scope of tear-film analysis and enable the estimation of key in vivo quantities, such as thinning rates, that typically require more complex modeling approaches (e.g., \cite{PengEtal2014,stapf2017duplex,DurschFLandThermal2017}).

In this article, we present results from fitting the evaporation-driven tear-film (TF) thinning model developed by Braun et al. \cite{braunTearFilm2018} and subsequently extended to two dimensions using proper orthogonal decomposition (POD) by Chen \cite{chenEvaporationdrivenTearFilm2024} to experimental fluorescence (FL) intensity data obtained from normal subjects. This analysis yields estimates of model parameters within tear-breakup regions that, to our knowledge, have not previously been determined in vivo. These results are expected to be of interest to both researchers and clinicians, advancing understanding of TF thinning and dry-spot formation and providing a useful reference for comparisons with data from dry eye patients.

The remainder of the paper is organized as follows. Section~\ref{sec:mathmodels} presents the mathematical model. Section~\ref{sec:methods} describes the inverse problem and numerical methods. Section~\ref{sec:exps} reports numerical results.

\section{Mathematical models}
\label{sec:mathmodels}
\subsection{Two-dimensional PDE}

Our focus is the model from Chen et al. \cite{chenEvaporationdrivenTearFilm2024}, which is a two-dimensional version of the models derived by Braun et al. \cite{braunTearFilm2018} and used
by Luke et al. \cite{lukeParameterEstimation2020} for local TBU dynamics (see \eqref{eq:A1}--\eqref{eq:AJ} in the appendix).

The TF is modeled as a Newtonian fluid over a flat corneal surface at the plane $z=0$. Because the film is thin, there is a separation of scales, and lubrication theory may be applied \cite{ODB97,CrasMat09review}.  The result of applying this perturbation approach is that an approximate velocity field is found, and the depth-averaged velocities over $0 < z < h(x,y,t)$ appear in a high order PDE for the TF thickness $h$.  The transport of solutes inside the tear film have been derived \cite{braunDynamicsTear2012,li2Dosmofluor} using an application of the theory developed by Jensen and Grotberg \cite{JenGrot93}. The key variables in the system are shown in \autoref{tab:variables}. These quantities are nondimensionalized according to
\begin{align}
x'&=\ell x , \quad y'=\ell y , \quad z'=dz, \quad t'=\frac{d}{v_\text{max}}t, \quad h'=dh,\quad u'=\frac{v_\text{max}}{\epsilon}u. \\
v'&=\frac{v_\text{max}}{\epsilon}v, \quad w'=v_\text{max}w, \quad J'=\rho v_\text{max}J, \quad c'=c_0 c,\quad f'=f_{\text{cr}}f,
\label{eq:TF2D_scale}
\end{align}
where primes denote dimensional quantities and the relevant physical parameters are given in \autoref{tab:physical}. The necessary nondimensional parameters are given in \autoref{tab:parameters}, where $\phi$ shows later in the FL intensity equation.

\begin{table}
    \centering
    \begin{tabular}{c c c}
     Parameter & Expression & Value  \\ \hline
     $\epsilon$ & $d/\ell$ & $8.3\times 10^{-3}$ \\[1mm]

     $P_{c}$ & $({P_{0} V_{w} c_{0}})/({v_\text{max}})$ & $0.392$  \\[1mm]
      $\text{Pe}_{f}$ & $({v_\text{max} \ell})/({\epsilon D_{f}})$ & $27.7$ \\[1mm]
      $\text{Pe}_{c}$ & $({v_\text{max} \ell})/({\epsilon D_{0}})$ & $6.76$ \\[1mm]
     $\phi$ & $\epsilon_{f} f_{\text{cr}} d$ & $0.417$
    \end{tabular}
    \caption{Typical values of nondimensional parameters that appear in the model \eqref{eq:pdeh}--\eqref{eq:pdef} using the parameters in \autoref{tab:physical}. The parameter $\phi$ appears in FL intensity equation \eqref{eq:FLI}.}
    \label{tab:parameters}
\end{table}

The resulting nondimensional system for $h,p,c,f$ is
\begin{align}
\partial_t h  + \partial_x (h\overline{u})+\partial_y (h\overline{v})&= -J + P_c(c-1), \label{eq:pdeh}\\
\overline{u} &=-\frac{h^2}{12}\partial_x p, \label{eq:pdeu}\\
\overline{v} &=-\frac{h^2}{12}\partial_y p, \label{eq:pdev}\\
p &=-\partial_x^2 h - \partial_y^2 h, \label{eq:pdep}\\
h(\partial_t c +\overline{u}\partial_x c + \overline{v}\partial_y c) &= \text{Pe}_c^{-1}(\partial_x (h\partial_x c)+\partial_y (h\partial_y c))+Jc-P_c(c-1)c, \label{eq:pdec}\\
h(\partial_t f +\overline{u}\partial_x f + \overline{v}\partial_y f) &= \text{Pe}_f^{-1}(\partial_x (h\partial_x f)+\partial_y (h\partial_y f))+Jf-P_c(c-1)f. \label{eq:pdef}
\end{align}
The variables $\overline{u}$ and $\overline{v}$ are depth-averaged transverse fluid velocities, and $p$ is the pressure inside the tear film.

\begin{table}[tbp]
    \centering
    \begin{tabular}{cl}
        Variable & Meaning \\ \hline
        $x$, $y$ & transverse spatial dimensions \\
        $z$ & depth dimension \\
        $t$ & time \\
        $h(x,y,t)$ & TF thickness \\
        $u(x,y,t)$, $v(x,y,t)$ & transverse fluid velocities \\
        $p(x,y,t)$ & pressure \\
        $J(x,y)$ & evaporation rate \\
        $c(x,y,t)$ & osmolarity \\
        $f(x,y,t)$ & fluorescein concentration
    \end{tabular}
    \caption{Variables in the two-dimensional model.}
    \label{tab:variables}
\end{table}

\begin{table}[tbp]
    \centering
    \begin{tabular}{c l c l}
      Parameter & Description & Value & Reference \\
     \hline
     $\mu$ & Viscosity & \qty{1.3e-3}{\Pa\sec} & Tiffany \cite{tiffanyViscosityHumanTears1991} \\
     $\sigma_{0}$ & Surface tension & \qty{0.045}{\N\per\meter} & Nagyov\'a and Tiffany \cite{nagyovaComponentsResponsibleSurface1999} \\
     $\rho$& Density & \qty{e3}{\kg\per \meter\tothe{3}} & Water \\
     $d$ & Initial TF thickness & \qty{4.5}{\micro\meter} & Calculated \\
     $\ell$ & $(\sigma_0/ \mu / v_{\max})^{1/4} d$ & \qty{0.54}{\milli \meter} & Calculated \\
     $v_{\text{max}}$ & Peak thinning rate & \qty{10}{\micro\meter\per\minute} & Nichols et al. \cite{nicholsThinningRatePrecorneal2005a} \\
     $V_{w}$ & Molar volume of water & \qty{1.8e-5}{\meter\cubed \per \mol} & Water \\
     $D_{f}$ & Diffusivity of fluorescein & \qty{0.39e-9}{\meter\squared \per \sec} & Casalini et al. \cite{casaliniDiffusionAggregationSodium2011} \\
     $D_{o}$& Diffusivity of salt & \qty{1.6e-9}{\meter\squared \per \sec} & Riquelme et al. \cite{riquelmeInterferometricMeasurementDiffusion2007} \\
     $c_{0}$& Isotonic osmolarity & \qty{300}{\Osm \per \meter\cubed} & Lemp et al. \cite{lempTearOsmolarityDiagnosis2011}\\
     $P_{0}$& Permeability of cornea & \qty{12.1}{\um\per\sec} & Braun et al. \cite{braunDynamicsFunctionTear2015}\\
     $\epsilon_{f}$& Napierian extinction coefficient & \qty{1.75e7}{\L \per \meter \per \mol} & Mota et al. \cite{motaSpectrophotometricAnalysisSodium1991}\\
     $f_{cr}$& Critical FL concentration  &  $0.2\%$ & Webber and Jones
     \cite{webberContinuousFluorophotometricMethod1986}
    \end{tabular}
    \caption{Physical parameters (dimensional) used in the governing equations.}
    \label{tab:physical}
\end{table}

We assume that the simulation takes place over a section of the cornea that is not close to the eyelids and limbus. Because we are not interested in the effects of these boundaries, we assume periodic spatial behavior on all the dependent variables. We also assume that the simulation begins after the eye opens and that all the dependent variables are initially uniform:
\begin{equation}
\label{eq:ic}
    h(x,y,0) = c(x,y,0) = 1,\quad f(x,y,0) = f_0, \quad p(x,y,0)=0,
\end{equation}
where $f_0$ is the FL concentration normalized to the critical concentration $f_{\text{cr}}$.

The evaporation rate $J$ is our primary input to the model, and it drives all the dynamics. Inhomogeneities in an in vivo lipid layer, sitting atop the aqueous layer of the tear film, are presumed to cause local increases in the evaporation rate, leading to local decreases in $h$ and corresponding increases in the solute concentrations. We represent the spatial variation of $J$ as one or more localized peaks:
\begin{equation} \label{eq:Jm}
J(x,y) = v_b + \sum_{k=1}^{K} (a_k-v_b)\,G\left( \frac{x-x_k}{x_{w,k}}, \frac{y-y_k}{y_{w,k}} \right),
\end{equation}
where $v_b = v_{\text{min}}/v_{\text{max}}$ is a baseline value, $(x_k,y_k)$ is the center of the $k$th peak, $a_k>v_b$ is the height of the $k$th peak, $x_{w,k}$ and $y_{w,k}$ are characteristic widths of peak $k$, and $G$ is the Gaussian
\begin{equation}
    \label{eq:gaussian}
    G(x,y) = \exp\left[-(x^2 + y^2) / 2 \right].
\end{equation}

The FL intensity $I$ is obtained via \cite{webberContinuousFluorophotometricMethod1986,braunModelTearFilm2014}
\begin{equation}
\label{eq:FLI}
I=I_0\frac{1-\exp(-\phi fh)}{1+f^2},
\end{equation}
where $I_0$ is a normalization coefficient and $\phi$ is the nondimensional Napierian extinction coefficient in \autoref{tab:parameters}.

Given an evaporation function $J(x,y)$, we solve the system \eqref{eq:pdeh}--\eqref{eq:pdef} to obtain $h$, $p$, $c$, and $f$ as functions of space and time. These can then be inserted into~\eqref{eq:FLI} to find the fluorescent intensity.

\section{Methods}
\label{sec:methods}

We address the inverse problem for the evaporation-driven tear film model. Specifically, given an FL intensity video in the vicinity of a likely TBU, we would like to estimate osmolarity $c$, TF thickness $h$, and FL concentration $I$ over space and time. We will achieve this by estimating the parameters in a simple evaporation function in the mathematical model via optimized matching to an observation of intensity.

We perform parameter estimation on experimental data taken in a study conducted at Indiana University \cite{awisi-gyauChangesCornealDetection2019a}. The study received approval from the Biomedical Institutional Review Board of Indiana University. Images were collected using fluorescein imaging from 25 participants. Each subject underwent a screening process before completing two separate visits, with ten imaging trials conducted at each visit. Subjects were excluded if they wore contact lenses or had been diagnosed with dry eye syndrome by a clinician.

A trial is defined as a sequence of eye images captured every $0.2$ or $0.25$ seconds, depending on the frame rate. At the beginning of each trial, a 2 microliter drop of 2\% sodium fluorescein solution was applied to the eye \cite{carlson2016clinical}. The eye was illuminated using a cobalt blue excitation filter (494 nm), while a Wratten no. 8 yellow barrier filter was positioned along the imaging axis. The fluorescein in the TF emitted green light (521 nm) that was recorded \cite{book}.

The video recording began as the subject blinked three times to evenly distribute the fluorescein across the tear film. During these initial blinks, the light source was set to a low intensity, and a custom MATLAB algorithm \cite{wuEffectsIncreasingOcular2015} was used to estimate the initial fluorescein concentration, which was assumed to be uniform across the cornea. After the third blink, the light intensity was increased to a predefined high setting, and the subject was instructed to keep their eye open for as long as possible. The trial ended with the subject's next blink. Each trial yields a movie: a sequence of images starting with the low-light blinks and ending with the final blink. For results presented in this paper, S, V and T refer to subject number, visit number and trial number respectively.

\subsection{Optimization problem}
\label{sec:evap_func}
The optimization problem is
\begin{equation}
\mathop{\mathrm{argmin}}\limits_{p}
\left\| I_{\text{th}}(x,y, t;\, p) - I_{\text{ex}}(x,y,t) \right\|_2^2,
\label{eq:opt_norm}
\end{equation}
where $\|\cdot\|_2$ denotes the discrete $\ell^2$ norm taken over all spatial grid points under the region we selected
and all sampled times. $I_{\text{th}}$ represents the theoretical FL intensity that is computed via our 2D model \eqref{eq:pdeh}--\eqref{eq:pdef} and \eqref{eq:Jm}, $I_{\text{ex}}$ represents the experimental FL intensity that is obtained from the FL imaging data, and $p$ represents a parameterization of the evaporation function $J(x,y)$. The relative error at a time $t_j$ is given by
\begin{equation}
\mathrm{RelErr}(t_j)
=
\frac{
\left\|
I_{\mathrm{th}}(x,y,t_j;p)
-
I_{\mathrm{ex}}(x,y,t_j)
\right\|_2
}{
\left\|
I_{\mathrm{ex}}(x,y,t_j)
\right\|_2
},
\label{eq:rel_err_eqn}
\end{equation}

\subsection{Elliptic spot representation}

While a TBU instance with circular symmetry can effectively be modeled by a 1-D PDE \cite{lukeParameterEstimation2020}, a spot that lacks circular symmetry may be better represented by an elliptical evaporation function. One possibility for parameterization of $J$ is to use \eqref{eq:Jm}, with the parameters $v_b, a_1, x_w, y_w, x_k,y_k$, plus a rotation angle $\theta$. We have found it more robust, however, to replace $(x_w, y_w, \theta)$ with the focal vector $F=(f_x, f_y)$, which is the vector from the center to a focus, and eccentricity $e$. The angle of the ellipse's major axis is $\theta=\tan^{-1}({f_y}/{f_x})$, and we have the semi-major axis $a=c/e$ and semi-minor axis $b=\sqrt{a^2-c^2}$, where $c=||F||$ is the distance from the ellipse center to each focus. The evaporation distribution for an ellipse is therefore
\[
J(x,y) = \beta \cdot v_b + (a_1-\beta \cdot v_b)\,\hat{G}\left(x-x_1, y-y_1 \right),
\]
where
\begin{gather}
\hat{G}(x,y) = \exp\left[-\frac{1}{2}Q(x,y) \right], \\
Q(x,y) = \begin{bmatrix} x & y  \end{bmatrix} A \begin{bmatrix} x  \\ y \end{bmatrix}, \\
A =
\begin{bmatrix}
\cos\theta & -\sin\theta \\
\sin\theta & \cos\theta
\end{bmatrix}
\begin{bmatrix}
\frac{1}{a^2} & 0 \\
0 & \frac{1}{b^2}
\end{bmatrix}
\begin{bmatrix}
\cos\theta & \sin\theta \\
-\sin\theta & \cos\theta
\end{bmatrix}.
\end{gather}
The optimization parameters for the ellipse are $p = (v_{b},\, a_{1},\, f_x,\,f_y,\,x_{1},\,y_{1},\,e,\,\beta)$.  We introduce the additional parameter $\beta$ because $v_b$ depends on other physical constants ${P}_{c}$, $\text{Pe}_{c}$, and $\text{Pe}_{f}$. To retain independent control over the background level, we therefore include a separate parameter $\beta$ to adjust the background values.

\subsection{Numerical methods}
To solve the system \eqref{eq:pdeh}--\eqref{eq:pdef}, we first note that equation \eqref{eq:pdef} describes how the evolution of FL concentration $f$ depends on $h$, $p$, and $c$, but those quantities do not in turn depend on $f$. Hence, given an evaporation function $J(x,y)$, we solve first the system \eqref{eq:pdeh}--\eqref{eq:pdec} to obtain $h$, $p$, and $c$, and then separately solve~\eqref{eq:pdef} to find $f$. We use the method of lines with a Fourier spectral collocation method~\cite{trefethenSpectralMethodsMatlab2000} in space on a uniform periodic grid on the domain $(-\pi, \pi]^2$. The number of grid points $m$ and $n$ in each dimension is chosen to be even. In this paper, we choose $m=n=40$. The resulting discretization of spatial terms in \eqref{eq:pdeh}--\eqref{eq:pdec} creates a differential--algebraic system (DAE) that is solved in Julia using the QNDF solver, an adaptive quasi-constant time step stiff method in the \textbf{DifferentialEquations} package \cite{rackauckas2017differentialequations} similar to backward differentiation formulas and using Shampine's accuracy-optimal kappa values as defaults \cite{shampine1997matlab}. A more efficient forward solver can be obtained using proper orthogonal decomposition (POD), which was implemented for this problem by Chen et al.~\cite{chenEvaporationdrivenTearFilm2024}. The central idea of POD is to accelerate large dynamical systems by projecting their solutions onto a low-dimensional subspace spanned by dominant modes, which are extracted from short-time solution snapshots using singular value decomposition (SVD).

For the optimization process, we use the POD method for the forward solver and a gradient-free optimization algorithm in the NLopt package \cite{NLopt}. In most cases, we obtained the best results from the PRAXIS algorithm, but in a few cases, Nelder-Mead was faster. On a $40\times40$ spatial grid, the typical computational time for one forward solution in the optimization process is about 15 seconds. We optimize for the ellipse parameters described in \ref{sec:evap_func}. Basic image processing of the final video frame allows us to make a good initial guess for the foci and eccentricity.

In practice, if $J$ is non-periodic, meaning the evaporation width might be too large, then it immediately returns a large penalty and the optimizer restarts. Because we are most interested in the center of image and because we modified the edges to create periodic data, we restrict the norm in \eqref{eq:opt_norm} to a rectangle in the center, for instance as shown on the right in \autoref{fig:s28_spot_I}. Details of data preprocessing can be found in Appendix \ref{secA3}.

The elapsed time during the captured FL video sets the timescale. Given a sequence of images where time $t$ goes from $0$ to $T'$ dimensionally, we use the scale in \eqref{eq:TF2D_scale} to obtain the nondimensional time as $0$ to $T$.

However, if our forward solver is not able to solve up to $t=T$ because the evaporation is too large, then the objective function returns a large penalty value so that the optimizer will choose other parameters.

\section{Numerical experiments}
\label{sec:exps}

In this section, we first present optimization results from synthetic data, then fit several cases with experimental data.

In order to test the feasibility of the optimization process, we start with results from synthetic data generated via \eqref{eq:pdeh}--\eqref{eq:pdef}. To generate data for the inverse problem, we solve on a $40\times40$ grid without use of POD acceleration.

\begin{figure}
    \centering
    \includegraphics[width=0.85\textwidth]{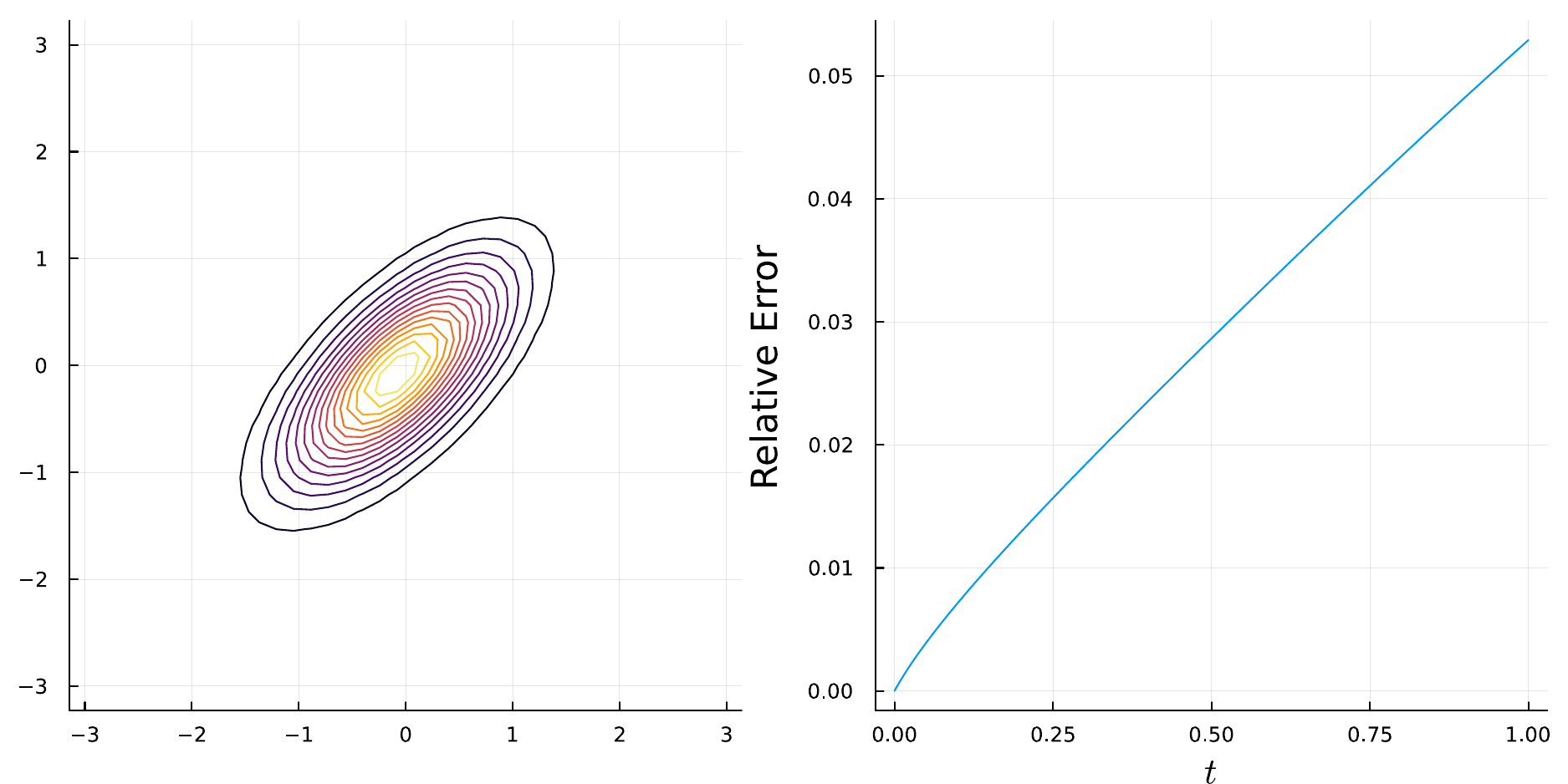}
    \caption{Left: contour plot of the evaporation function for an elliptic spot. Right: relative error of FL intensity for a 1-D radial fit lifted to 2-D.}
    \label{fig:synthetic_opt}
\end{figure}

The synthetic example is generated for an evaporation function with parameters $\beta=0.5, v_b=0.07, a_1=0.8, f_x=f_y=0.5, x_k=y_k=0, e=0.9$, as illustrated in \autoref{fig:synthetic_opt}. The 2-D optimization was initialized with $\beta=1, v_b=0.2, a_1=1.5, f_x=f_y=0.3, x_k=y_k=0.1, e=0.7$.  The optimization recovered the true parameter values up to $8$ digits in $83$ iterations with PRAXIS method.

We compare our full 2-D elliptical optimizations to the best results obtained in the same inverse problem by the 1-D model with circular symmetry, \eqref{eq:A1}--\eqref{eq:A5}, applied to radially averaged 2-D data, to find the parameters $v_{b}, r_{w}, a_{1},$ and $\beta$ in~\eqref{eq:AJ}. After $51$ iterations with the PRAXIS method, the optimized parameters are $v_b=0.06, r_w=0.7, a_1 = 0.54, \beta=0.8 $. \autoref{fig:synthetic_opt} shows that the error of the optimized 1-D radial solution grows to over 5\% by the end of the simulation.

\begin{figure}
    \centering
    \includegraphics[width=0.85\textwidth]{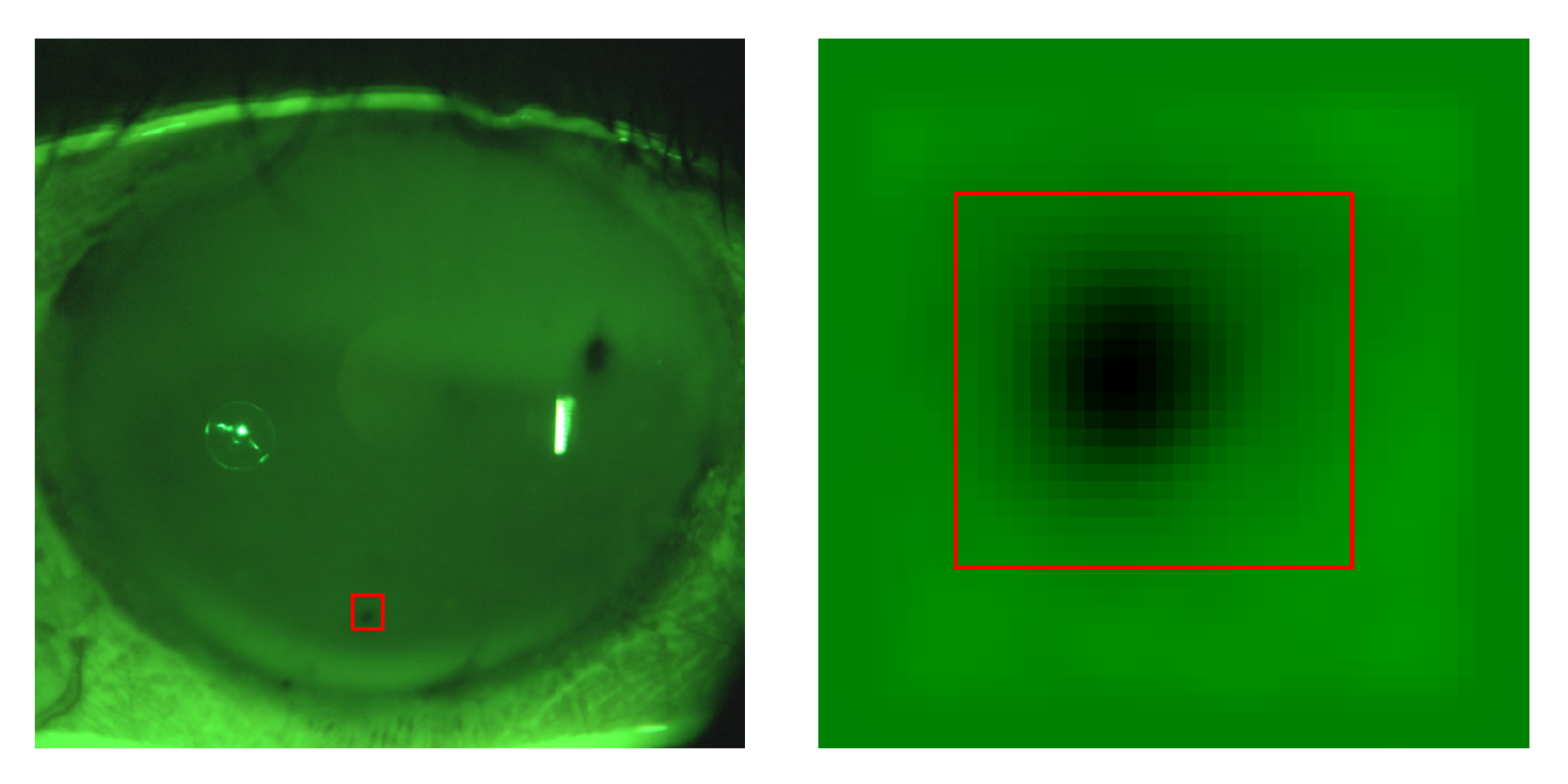}
    \caption{Case 1: Left: FL image of S28V1T3 with the local TBU region highlighted; Right: Smoothed data with the highlighted region for optimization.}
    \label{fig:s28_spot_I}
\end{figure}

Now we turn to the experimental data. We start with a nearly circular TBU case from the experimental data S28V1T3. In \autoref{fig:s28_spot_I} we show the global FL intensity image at the final time and highlight a likely TBU. We will denote it as Case 1. The optimized parameters are shown in \autoref{tab:opt_cases}, with a mildly eccentricity $e=0.19$. The optimized parameters using the 1D model are $v_b = 0.02$, $r_w=0.18$, $a_1=0.16$, and $\beta=2.02$.

\begin{figure}
    \centering
    \includegraphics[width=0.85\textwidth]{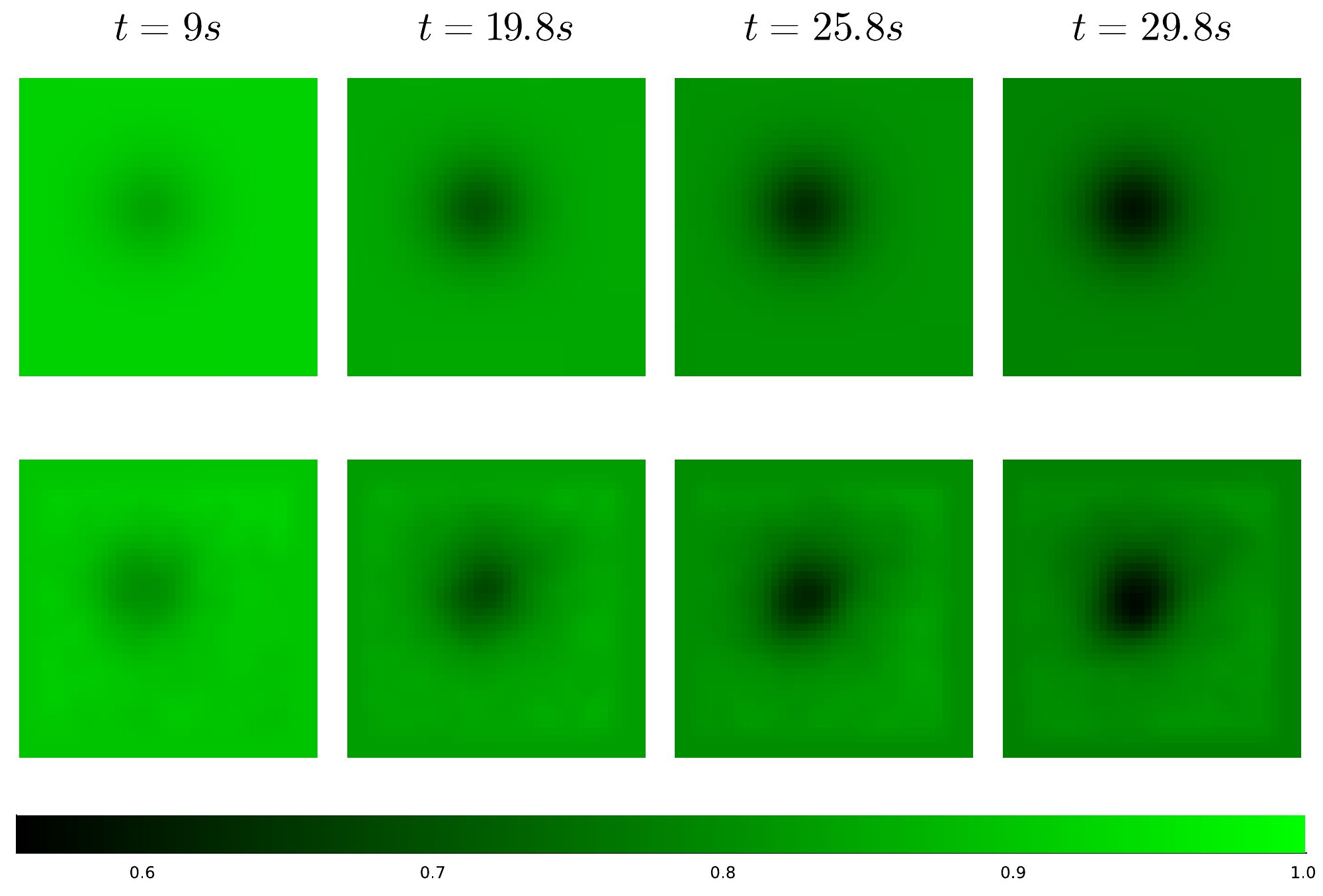}
    \caption{Top row: FL image of the fit with optimized parameters using an ellipse; Second row: Smoothed FL image of local elliptical spot in Case 1.}
    \label{fig:s28_spot_fit_I}
\end{figure}

\begin{figure}
    \centering
    \includegraphics[width=0.85\textwidth]{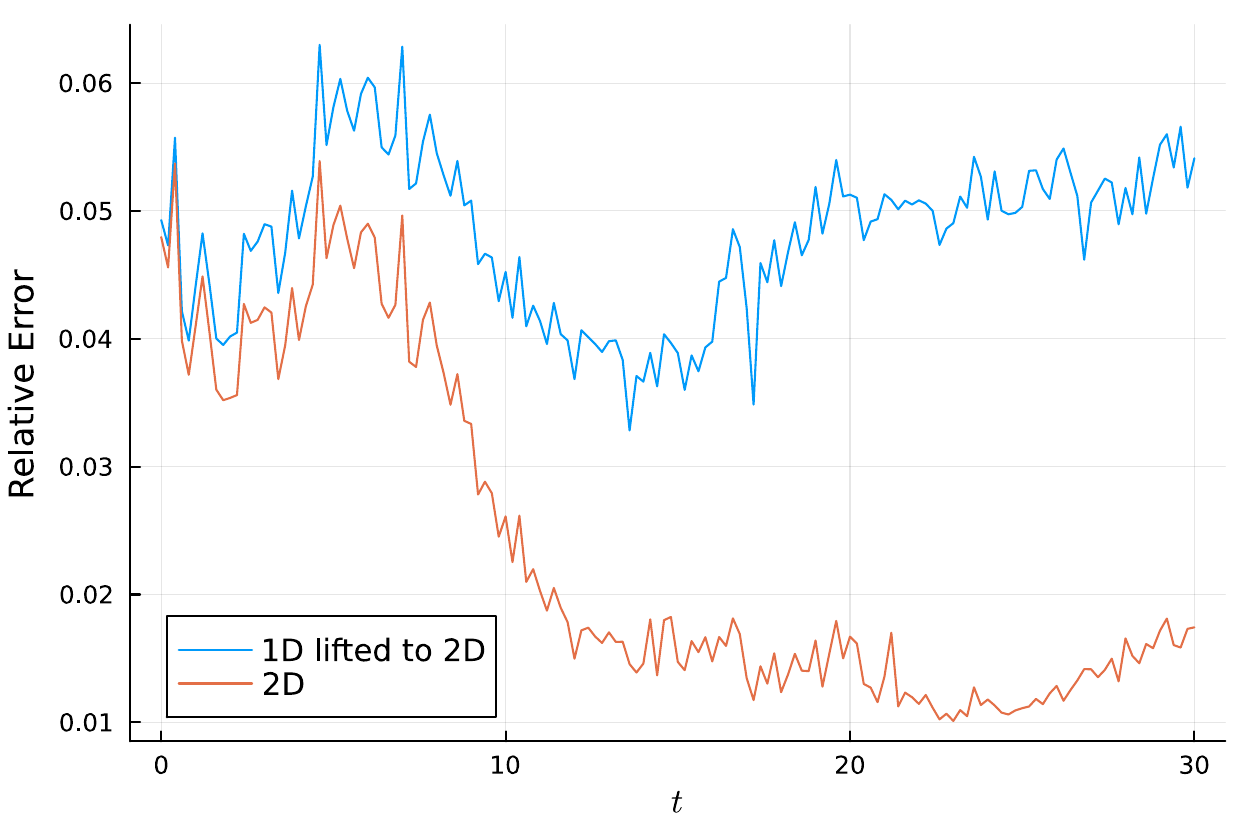}
    \caption{Red curve: Relative error of FL intensity for the 2D fit. Blue curve: Relative error of FL intensity for the 1D fit lifted to 2D. Case Number: Case 1.}
    \label{fig:s28_spot_err}
\end{figure}

\autoref{fig:s28_spot_fit_I} compares snapshots from the solution produced by the 2D fit with snapshots from the same times in the data. \autoref{fig:s28_spot_err} shows that the 2D fit is very good, with the relative residual below $2\%$ for most of the simulation time. The relative residual for the 1D fit is around $5\%$, which is not unreasonable.

\begin{figure}
    \centering
    \includegraphics[width=0.85\textwidth]{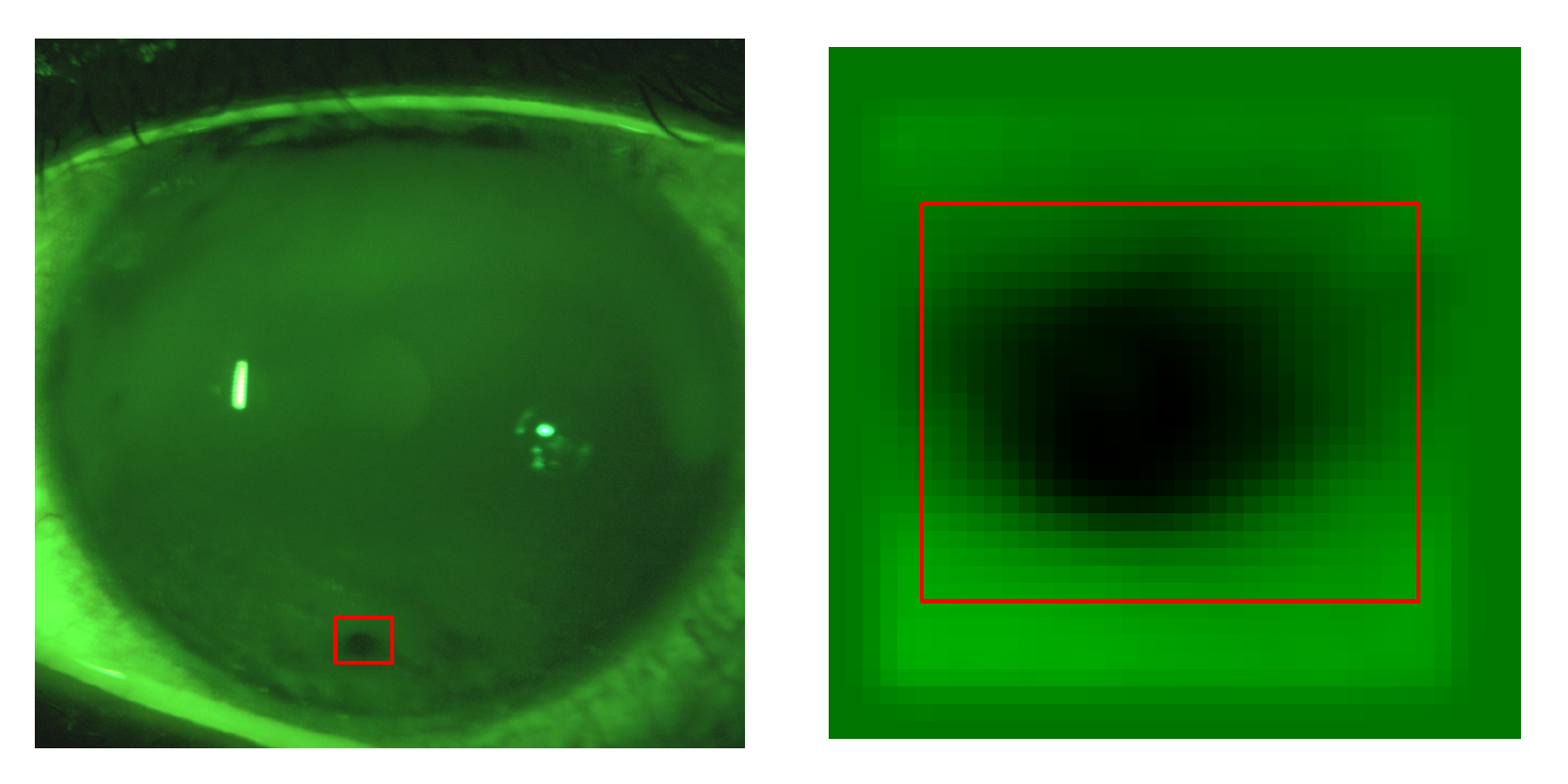}
    \caption{Case 2: Left: FL image of S27V1T8 with the local TBU region highlighted; Right: Smoothed data with the highlighted region for optimization}
    \label{fig:S27_plot}
\end{figure}

\begin{figure}
    \centering
    \includegraphics[width=0.85\textwidth]{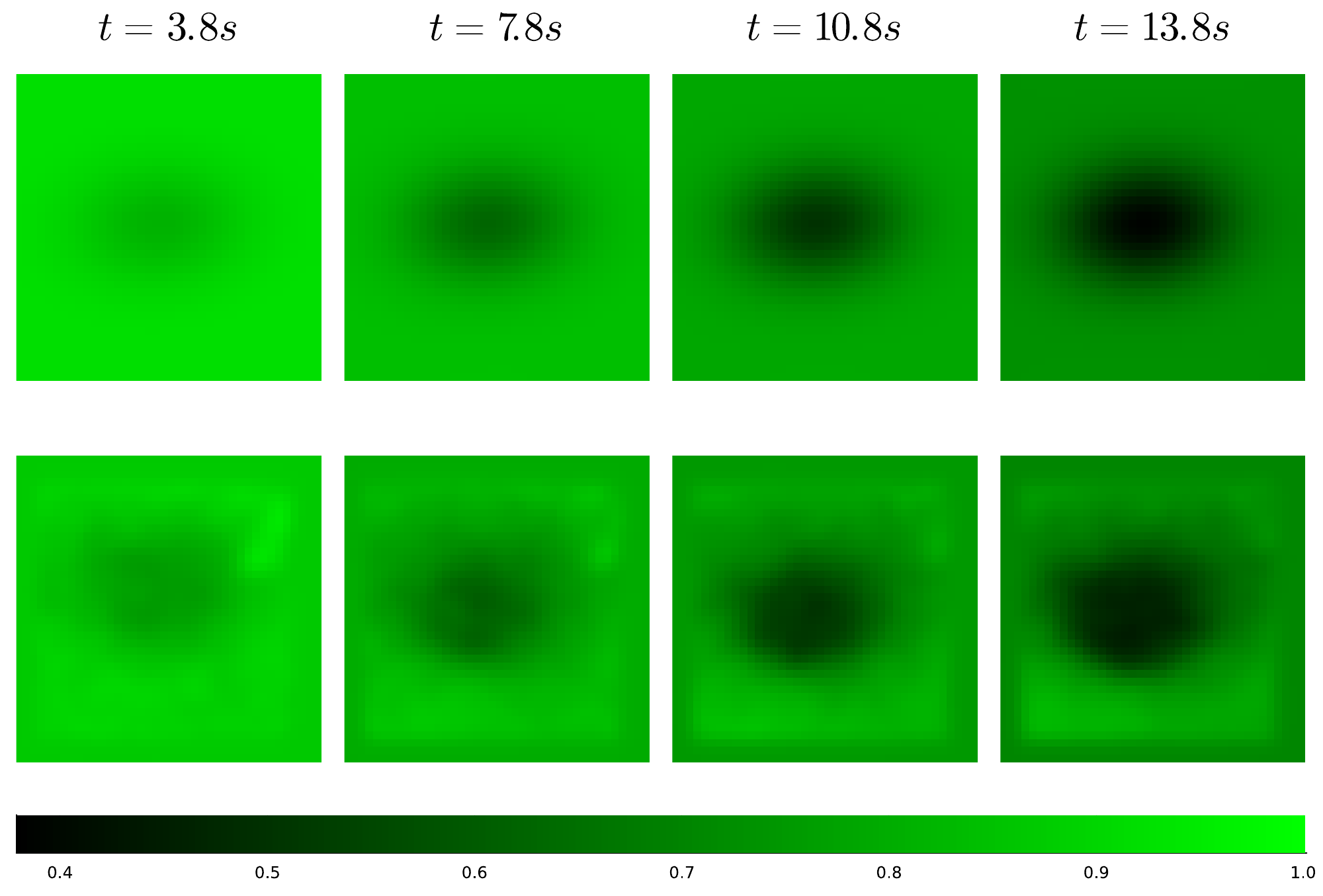}
    \caption{Top row: FL image of the fit with optimized parameters using an ellipse; Second row: Smoothed FL image of local elliptical spot in Case 2}
    \label{fig:S27v1t8b0_I}
\end{figure}

\begin{figure}
    \centering
    \includegraphics[width=0.85\textwidth]{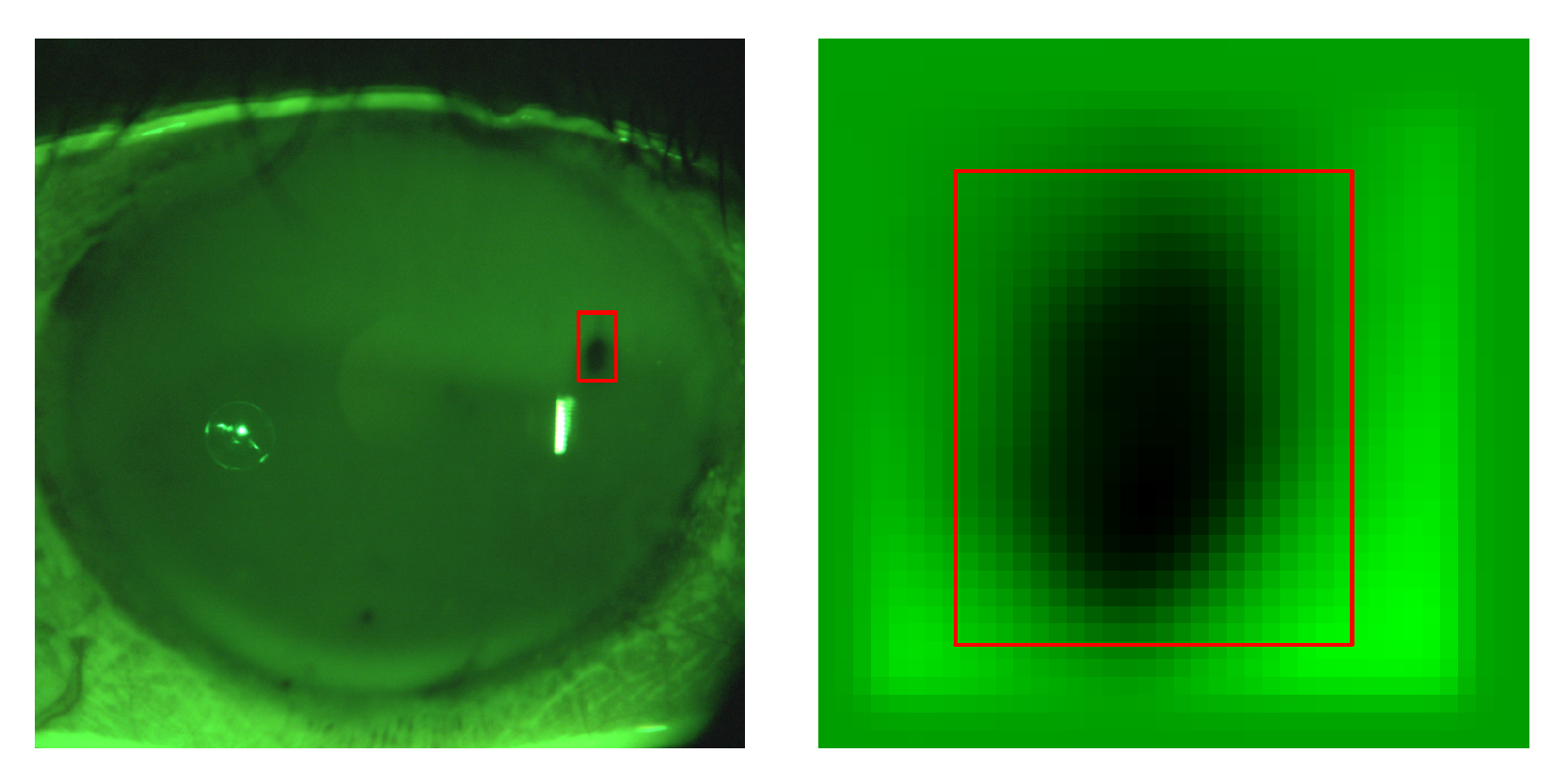}
    \caption{Case 3: Left: FL image of S28V1T3 with the local TBU region highlighted; Right: Smoothed data with the highlighted region for optimization}
    \label{fig:S28_ellp_I}
\end{figure}

\begin{figure}
    \centering
    \includegraphics[width=0.85\textwidth]{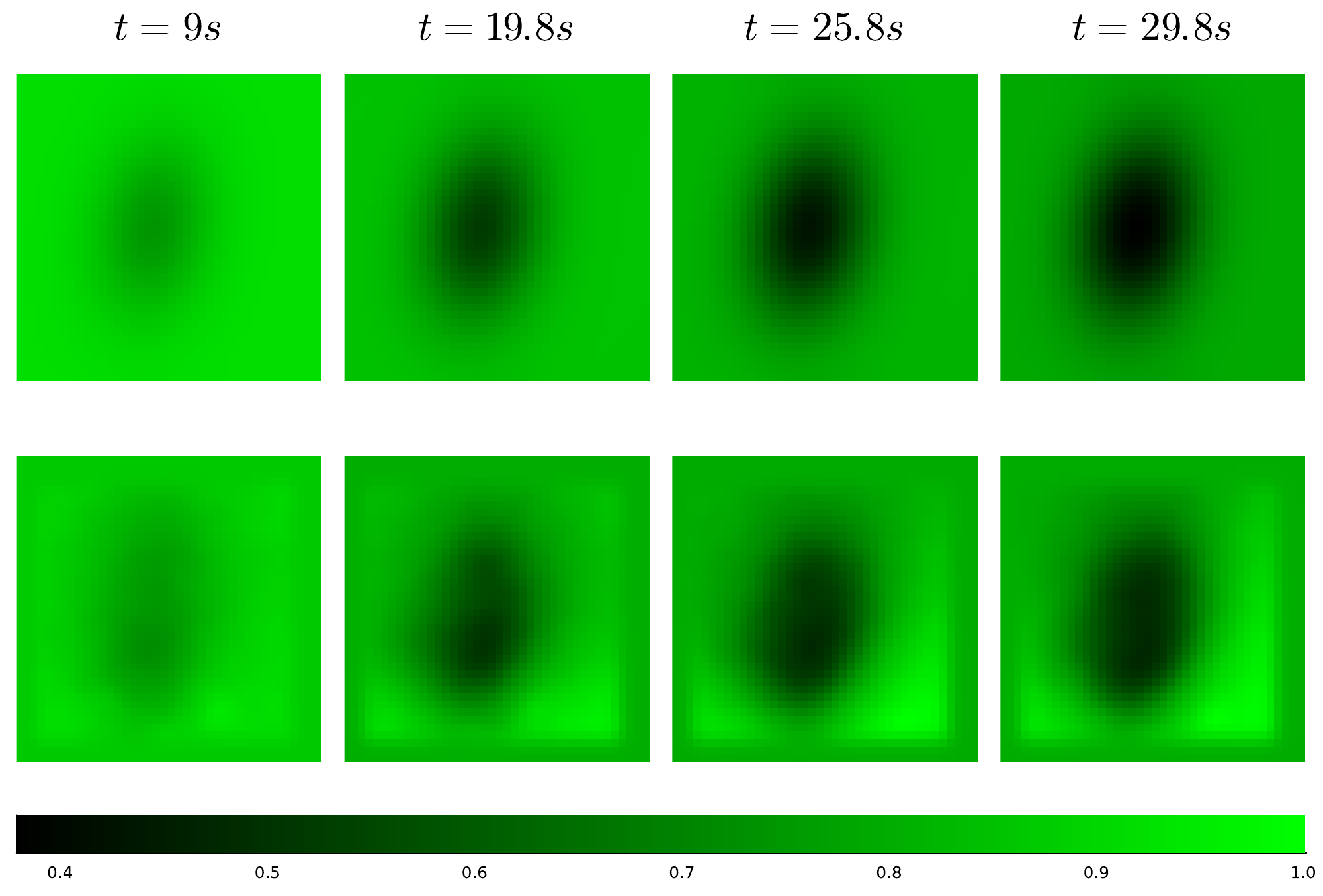}
    \caption{Top row: FL image of the fit with optimized parameters using an ellipse; Second row: Smoothed FL image of local elliptical spot in Case 3}
    \label{fig:S28_ellp}
\end{figure}

\begin{figure}[htbp]
  \centering
  \includegraphics[width=0.9\textwidth]{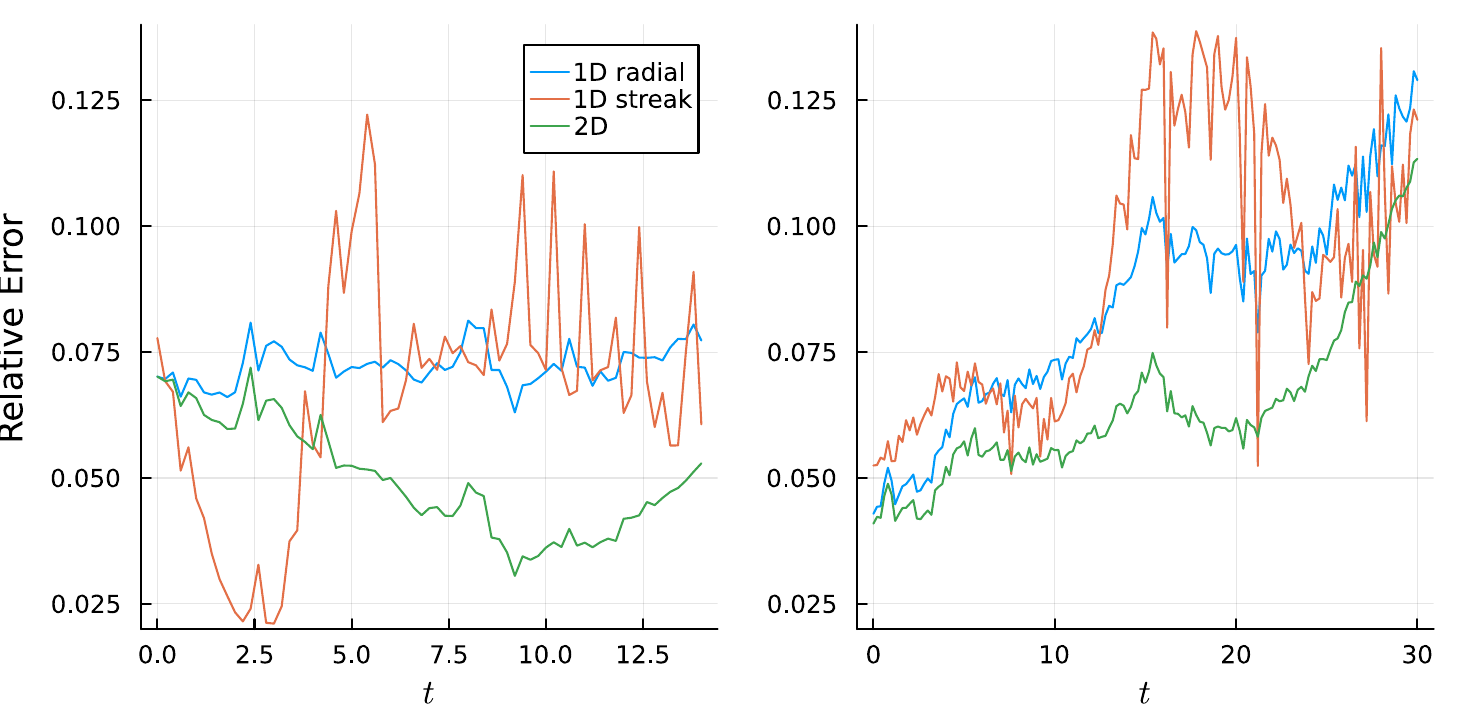}
  \caption{Left: Case 2; Right: Case 3; Blue and Red curves represent relative error of 1D radial and streak fit lifted to 2D, respectively.}
   \label{fig:s2728_rel_error}
\end{figure}

\autoref{fig:S27_plot} and \autoref{fig:S28_ellp_I} show final FL images from two other experimental videos and highlight TBU instances that are clearly more eccentric. We denote them as Case 2 and Case 3 respectively. As shown in \autoref{tab:opt_cases}, the eccentricity found in the best 2D fits were $e=0.75$ and $e=0.77$, respectively, quantifying the departure from circular symmetry. \autoref{fig:S27v1t8b0_I} and \autoref{fig:S28_ellp} show snapshots from the optimal fits to these spots, while \autoref{fig:s2728_rel_error} shows the relative residuals over time compared to the 1D fits. The optimized parameters using 1D radial model for these two cases are $(0.02,1.09,0.25,5.98)$ and $(0.48,1.26,2.98,1.73)$, respectively, for $(v_{b},\, r_{w},\, a_1,\,\beta)$. The left figure shows that the relative error is below $5\%$ most of the time for the 2D fit of Case 2, while the relative error for the 1D fit is below $8\%$. The right figure shows Case 3, where the relative error for the 1D fit exceeds $12\%$ at the final time, while the 2D fit remains at $10\%$.

\begin{table}[htbp]
\centering
\small
\setlength{\tabcolsep}{3pt}
\begin{tabular}{|c|c|c||c|c|c|c|c|c|}
\hline
Case & Iterations, Alg. & Type
& $f_x, f_y$
& $x_1, y_1$
& $v_b$
& $a_1$
& $e$
& $\beta$ \\
\hline

\multirow{2}{*}{Case 1}
& \multirow{2}{*}{105, PRAXIS}
& Initial
& (0.1, 0.1)
& (0, 0)
& 0.05
& 0.5
& 0.2
& 1 \\
& & Optimized
& (0.1, 0.13)
& (-0.26, 0.54)
& 0.08
& 0.37
& 0.19
& 2.2 \\
\hline

\multirow{2}{*}{Case 2}
& \multirow{2}{*}{127, PRAXIS}
& Initial
& (0.5, 0.1)
& (0, 0)
& 0.05
& 0.5
& 0.5
& 1 \\
& & Optimized
& (0.91, 0.05)
& (-0.07, 0.15)
& 0.014
& 0.28
& 0.75
& 7.68 \\
\hline

\multirow{2}{*}{Case 3}
& \multirow{2}{*}{144, PRAXIS}
& Initial
& (0.1, 1.5)
& (0, 0)
& 0.1
& 1
& 0.8
& 1 \\
& & Optimized
& (0.17, 1)
& (0.07, -0.22)
& 0.20
& 1.61
& 0.77
& 2.46 \\
\hline

\multirow{2}{*}{Case 4}
& \multirow{2}{*}{115, Nelder-Mead}
& Initial
& (0.1, 1)
& (0, 0)
& 0.02
& 1
& 0.7
& 1 \\
& & Optimized
& (0.22, 0.18)
& (0.07, -0.03)
& 0.025
& 0.42
& 0.14
& 8.9 \\
\hline

\end{tabular}
\caption{Initial guesses and optimized parameters for each case.}
\label{tab:opt_cases}
\end{table}

Since both Case 2 and Case 3 show high eccentricity, we can use 1D streak model \eqref{eq:B1}--\eqref{eq:B5} to do the fit as well. For Case 2, we extract the data horizontally along the minimum FL intensity over time, and for Case 3, we extract the data vertically along the minimum FL intensity. The red curve in \autoref{fig:s2728_rel_error} shows that the relative error of the 1D streak fit. For Case 2 it reaches the maximum $12.5\%$ while for Case 3 it is close to $14\%$. The error is small for Case 2 early in the simulation, possibly due to the uniformity at the early stages.

\begin{figure}
    \centering
    \includegraphics[width=0.85\textwidth]{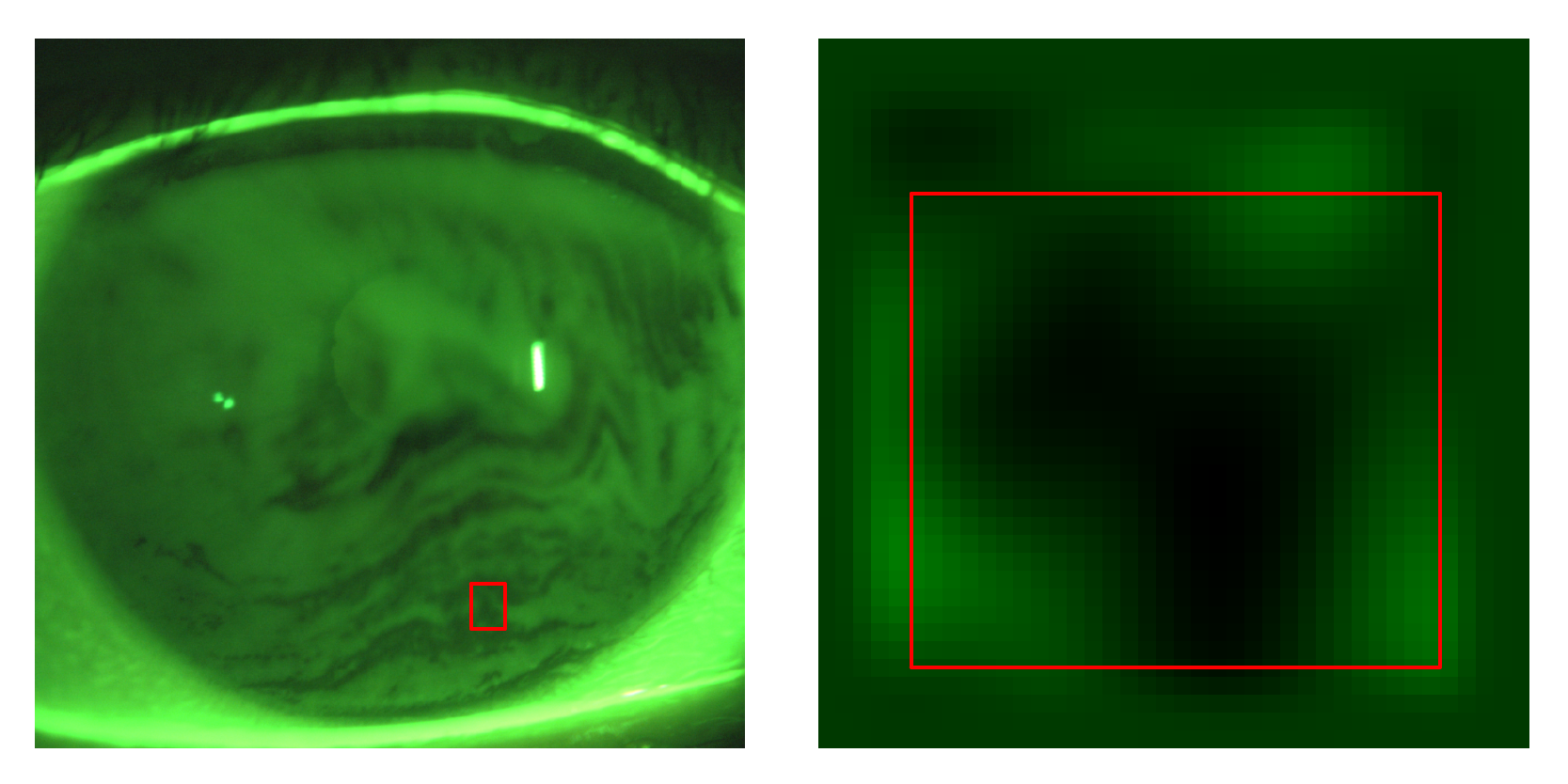}
    \caption{Left: FL image of Case 4 with the local TBU region highlighted; Right: Smoothed data with the highlighted region for optimization}
    \label{fig:S21_ellp_I}
\end{figure}

\begin{figure}
    \centering
    \includegraphics[width=0.85\textwidth]{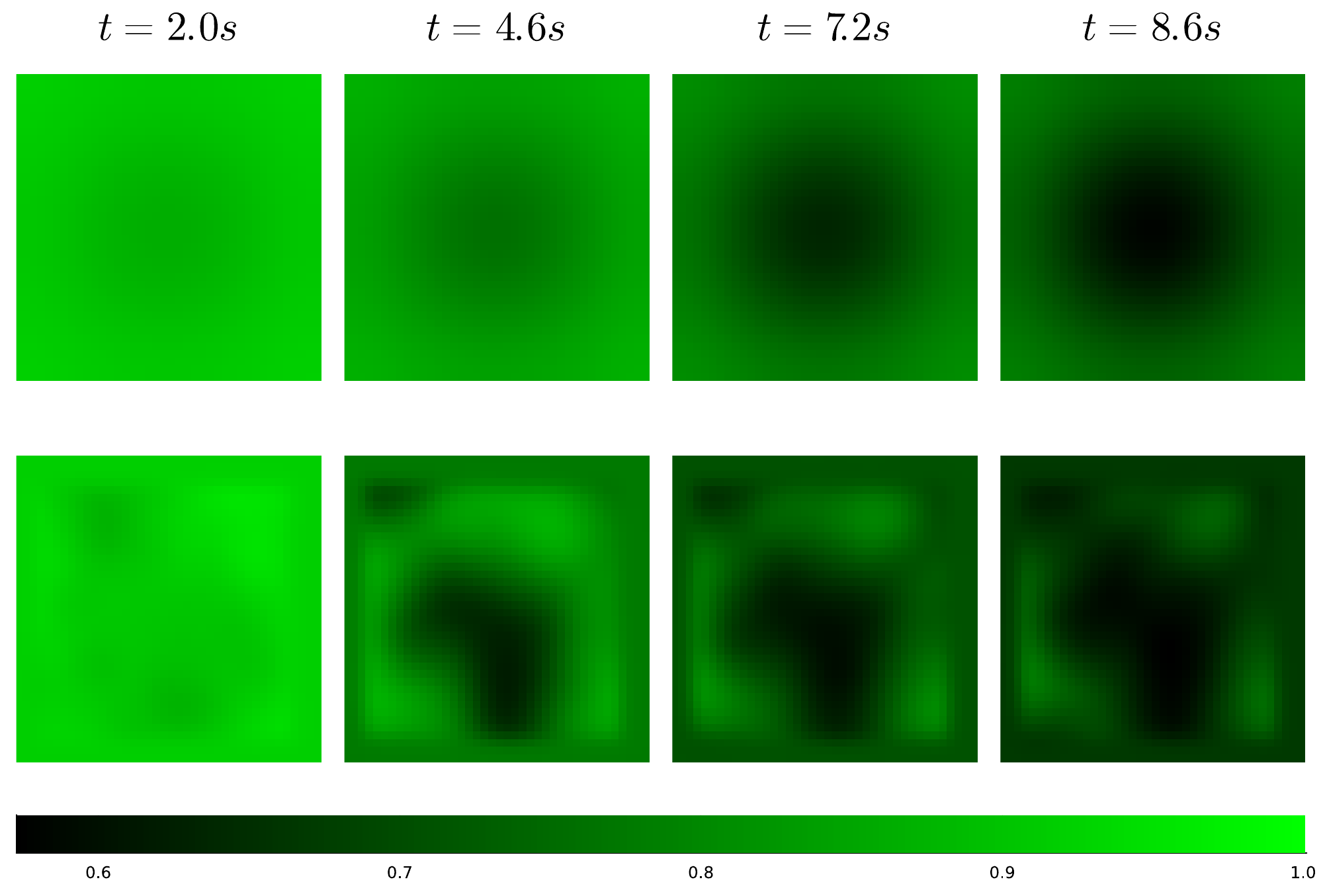}
    \caption{Top row: FL image of the fit with optimized parameters using an ellipse; Second row: Smoothed FL image of local elliptical spot in Case 4.}
    \label{fig:S21_I}
\end{figure}

\begin{figure}
    \centering
    \includegraphics[width=0.85\textwidth]{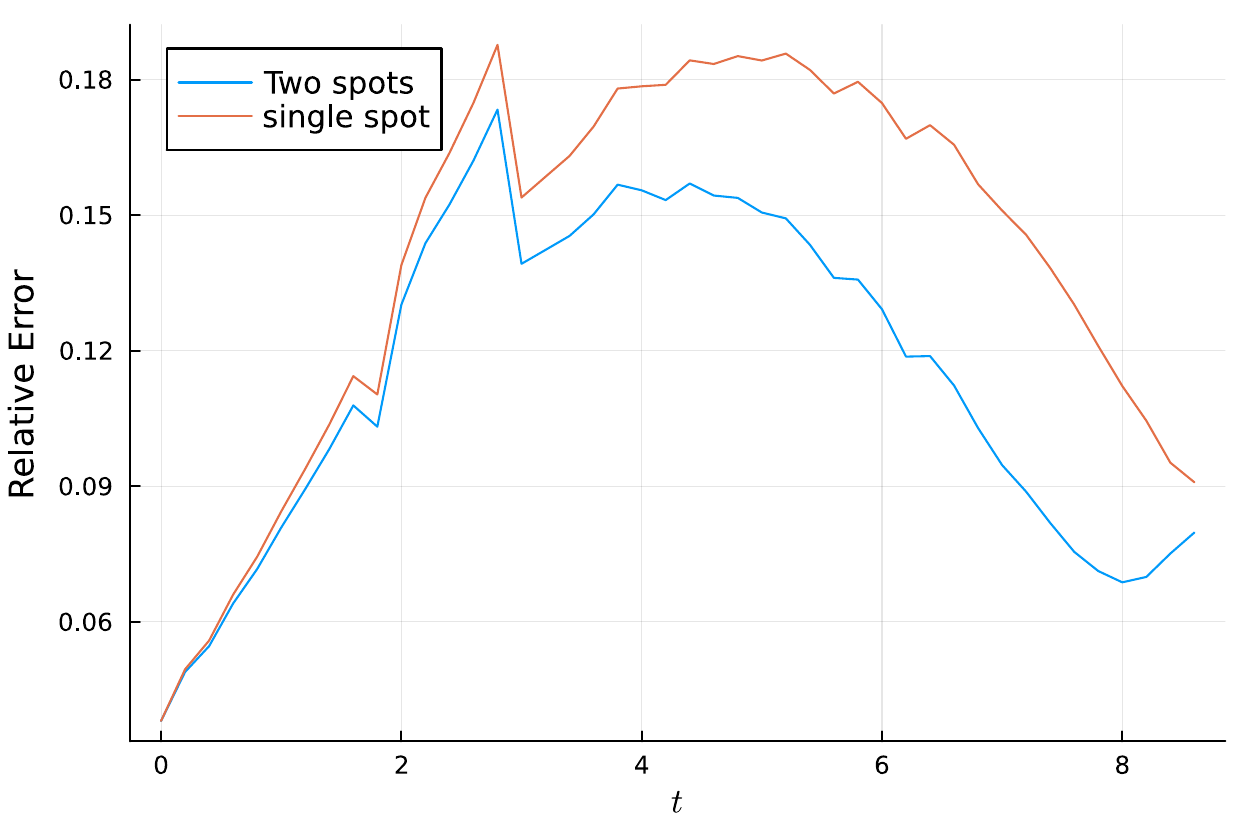}
    \caption{Red curve: Relative error with single spot fit for Case 4; Blue curve: relative error with two spots fit for Case 4.}
    \label{fig:S21_rel_error}
\end{figure}

\autoref{fig:S21_ellp_I} shows a severely nonuniform local TBU region in S21V1T1, which is denoted as our Case 4. As \autoref{fig:S21_I} shows, the optimization result is unsatisfactory using a single ellipse since the TBU shape is neither circular nor elliptic. We also use two elliptic spots for the fit, though the FL intensity image is not presented here since it is not satisfactory either. \autoref{fig:S21_rel_error} shows that using two spots for the fit has some improvements, specifically towards the end of the time. \autoref{tab:opt_cases} shows that the optimized eccentricity $e = 0.14$ which is close to a circular spot. This suggests that an alternative evaporation function which is more flexible should be implemented for such cases.

\begin{figure}
    \centering
    \includegraphics[width=0.85\textwidth]{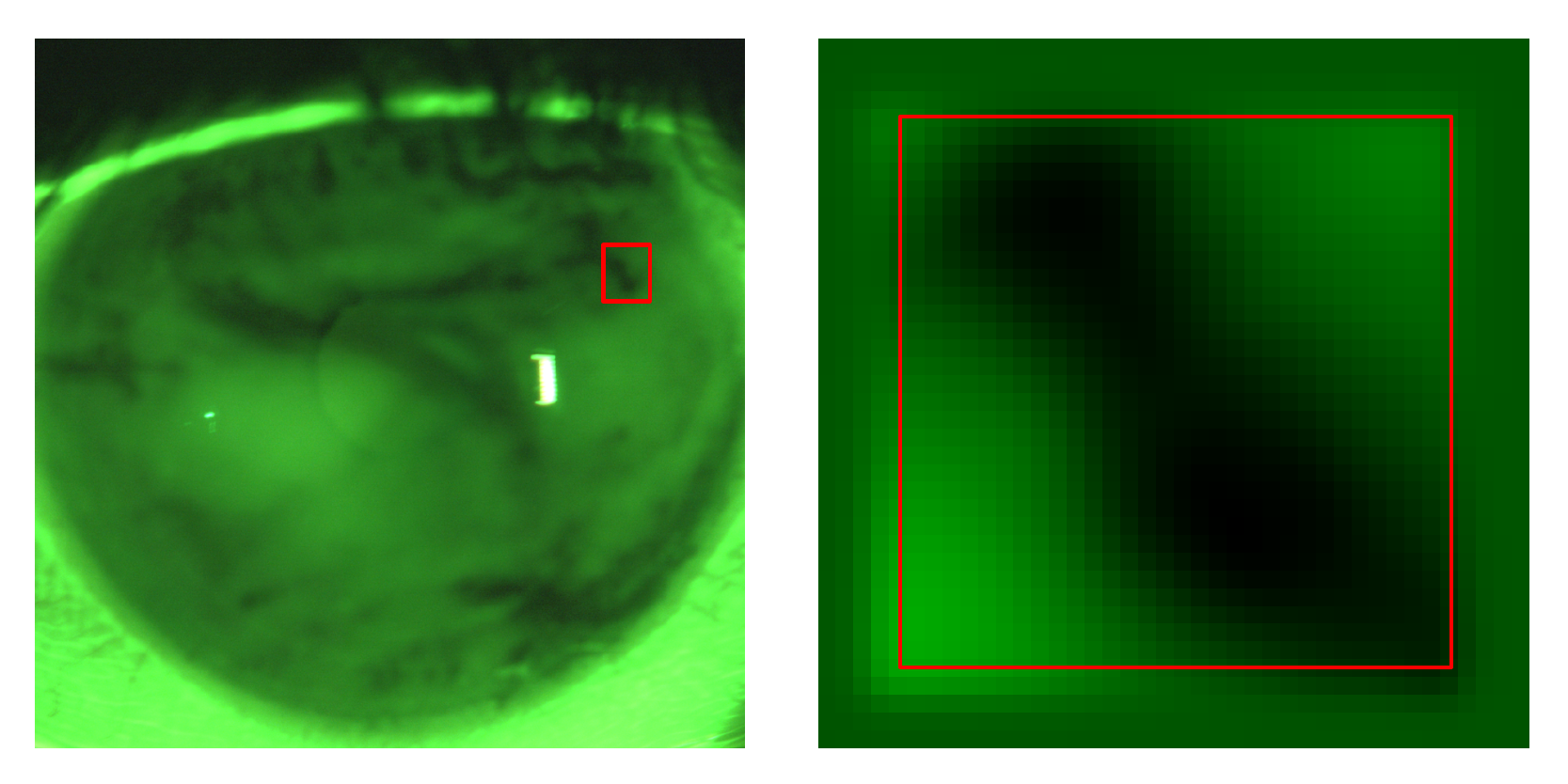}
    \caption{Case 5: Left: FL image of S1V1T5 with the local TBU region highlighted; Right: Smoothed data with the highlighted region for optimization}
    \label{fig:S1v1t5_ellp_I}
\end{figure}

\autoref{fig:S1v1t5_ellp_I} shows a local TBU region with multiple spots from S1V1T5. We denote it as Case 5. We use two elliptical representations of the evaporation function \eqref{eq:Jm}, and there would be $15$ parameters in total, namely $(v_{b},\, a_{1},\, a_2, \, f_{x_1},\,f_{y_1},\, f_{x_2}, f_{y_2}, x_{1},\,y_{1},\, x_{2}, \, y_{2},\, e_1,\, e_2, \, \beta_1,\, \beta_2)$. The initial guess we use is $(v_{b}=0.05,\, a_{1}=0.6,\, a_2=0.5, \, f_{x_1}=0.2,\,f_{y_1}=0.2,\, f_{x_2}=0.5, f_{y_2}=-0.5, x_{1}=-0.7,\,y_{1}=1.5,\, x_{2}=0.5, \, y_{2}=-0.5,\, e_1=0.5,\, e_2=0.7, \, \beta_1=1,\, \beta_2=2)$. The PRAXIS method does not converge for this case, and we use Nelder-Mead instead. It takes $224$ iterations in total. The optimized values are
\[
\begin{aligned}
v_{b}=0.048,\, a_{1}=1.06,\, a_{2}=0.56,\,
f_{x_1}=0.17,\, f_{y_1}=0.49,\,
f_{x_2}=0.48,\, f_{y_2}=-0.8, \\
x_{1}=-1.23,\, y_{1}=1.03,\,
x_{2}=0.55,\, y_{2}=-0.38,\,
e_{1}=0.51,\, e_{2}=0.88,\,
\beta_{1}=5.02,\, \beta_{2}=4.21
\end{aligned}
\]

\begin{figure}
    \centering
    \includegraphics[width=0.85\textwidth]{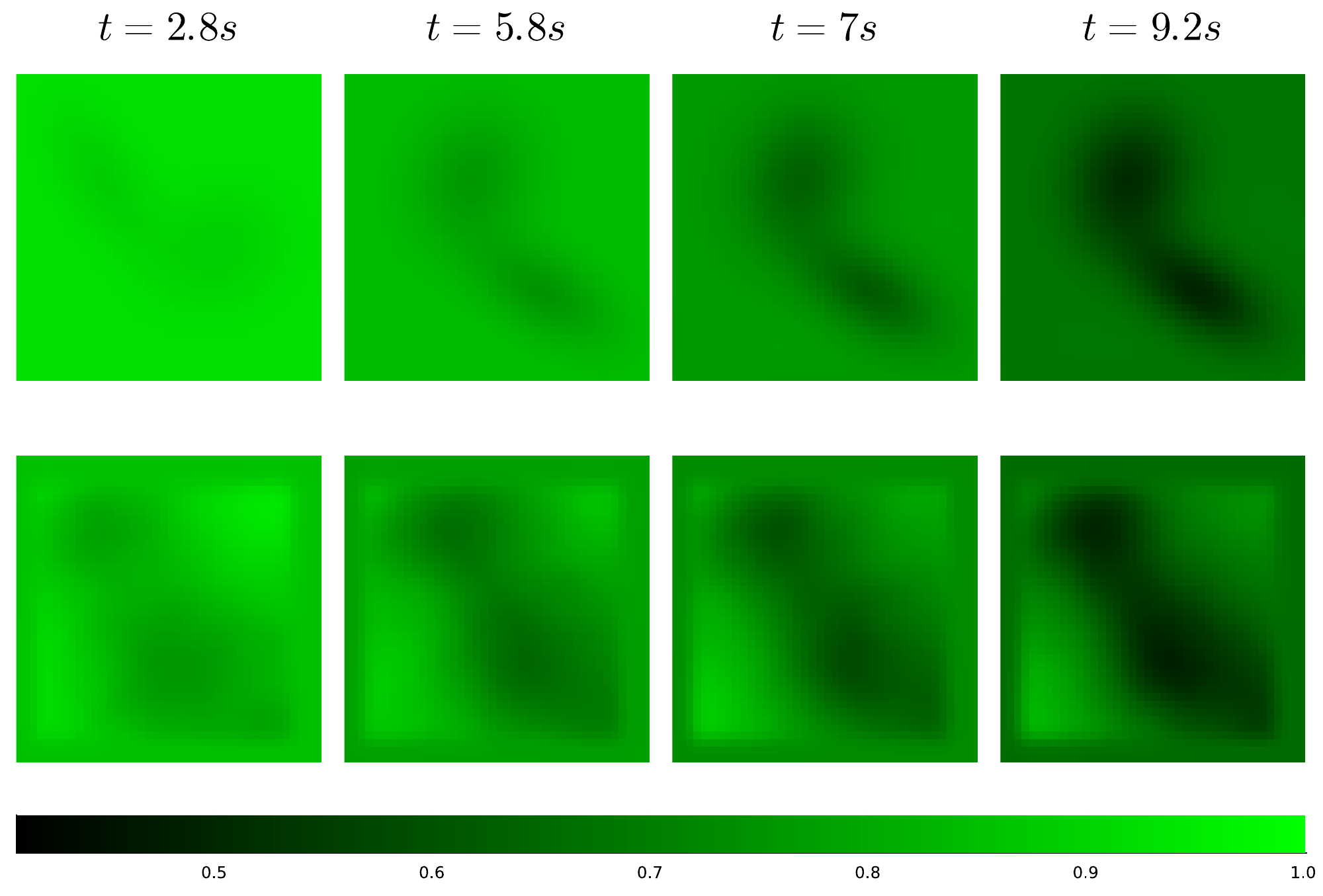}
    \caption{Top row: FL image of the fit with optimized parameters using an ellipse; Second row: Smoothed FL image of local elliptical spot in Case 5}
    \label{fig:S1v1t5_ellp}
\end{figure}

\begin{figure}[htbp]
  \centering
  \includegraphics[width=0.8\textwidth]{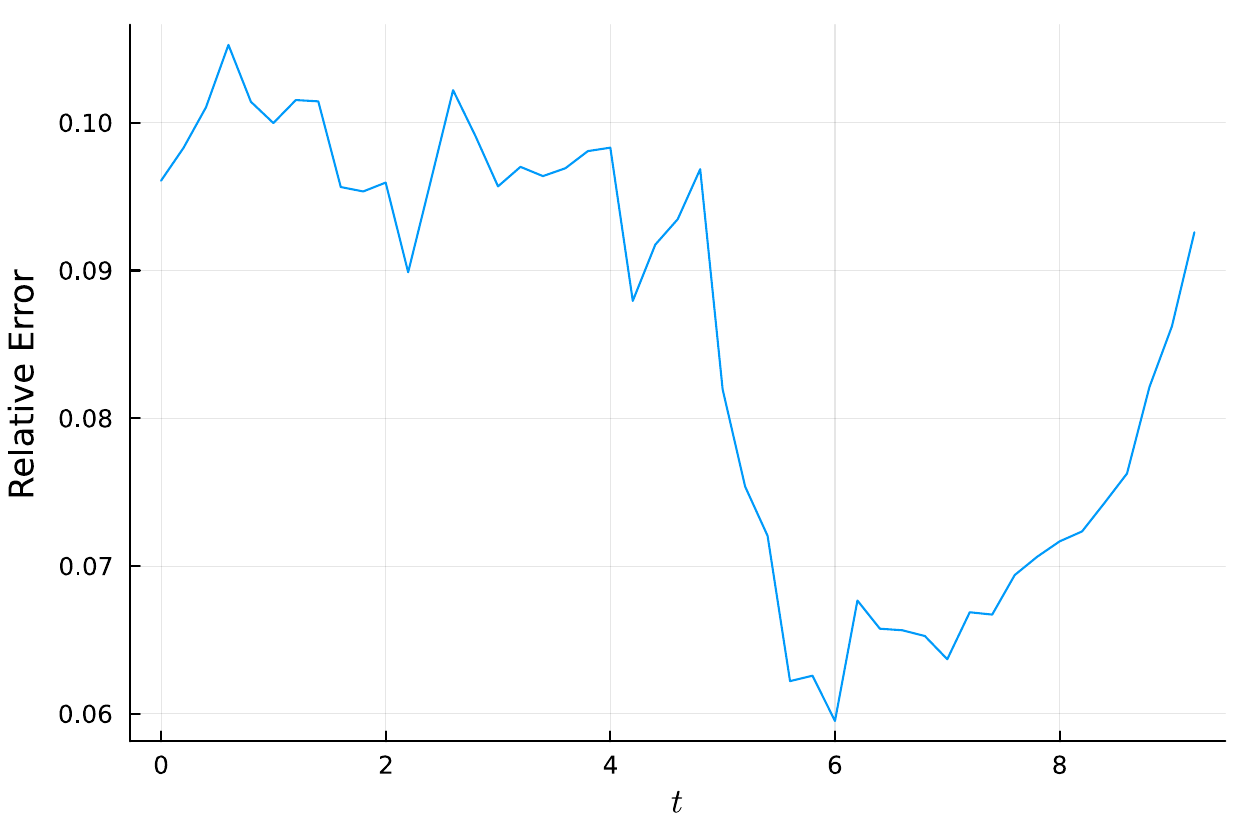}
  \caption{Relative error of Case 5.}
   \label{fig:s1v1t5_rel_error}
\end{figure}

\autoref{fig:S1v1t5_ellp} shows snapshots from the optimal fits to the data. The shapes are close, although there is still room for improvement. \autoref{fig:s1v1t5_rel_error} shows that the relative error is below $10\%$ most of the time.

\section{Discussion and future work}

The results presented here demonstrate that our two-dimensional PDE model, combined with POD reduction, provides a practical framework for estimating in vivo parameters from FL imaging data. Compared with one-dimensional streak or radial models, the 2D formulation of the model more faithfully captures asymmetric and multi-spot local TBU regions. Some cases are not well fit, likely because our evaporation distribution function is fixed. Future work will focus on extending the evaporation parameterization to allow more flexible spatial structure, incorporating flow-driven or lipid-dynamics effects. Application of this framework to dry-eye patient populations will be particularly important for determining whether the estimated evaporation and thinning parameters can serve as clinically meaningful markers of tear-film instability.

\begin{appendices}
\section{Axisymmetric model}\label{secA1}
For the circular case, we use the axisymmetric coordinates $(r',z')$ to denote the position and $u'=(u',w')$ to denote the fluid velocity. The tear film is modeled as an incompressible Newtonian fluid on $0 < r' < R_0$ and $0<z'<h'(r',t')$. The scalings and non dimensional parameters are similar in \autoref{tab:physical} and \autoref{tab:parameters}. The only difference is that $r'=\ell r$. Braun et al. derived the system of equations on the domain $0<r<R_0$ \cite{braunTearFilm2018,braunDynamicsFunctionTear2015}.

\begin{align}
\partial_t h  + \partial_r (rh\overline{u}) &= -J + P_c(c-1), \label{eq:A1}\\
\overline{u} &=-\frac{h^2}{12}\partial_r p, \label{eq:A2}\\
p &=-\frac{1}{r}\partial_r(r\partial_r h), \label{eq:A3}\\
h(\partial_t c +\overline{u}\partial_r c ) &= \text{Pe}_c^{-1}\frac{1}{r}\partial_r (rh\partial_r c)+Jc-P_c(c-1)c, \label{eq:A4}\\
h(\partial_t f +\overline{u}\partial_r f ) &= \text{Pe}_f^{-1}\frac{1}{r}\partial_r (rh\partial_r f)+Jf-P_c(c-1)f, \label{eq:A5}
\end{align}
The evaporation function $J$ is

\begin{equation} \label{eq:AJ}
J(r) = \beta\cdot v_b + (a-v_b) \exp\left[-(r/r_w)^2 / 2 \right].
\end{equation}

where $v_b$ is the ratio of $v_{\text{min}}$ over $v_{\text{max}}$, $r_w$ is the radius, and $a > v_b$ is the height of the peak. $\beta$ is a scaling constant.

\section{Full streak model}\label{secA2}
The linear case model is solved on the Cartesian coordinates $-\pi < x < \pi$ and $ 0 < z < h(x,t)$. More details about derivation can be found in \cite{braunTearFilm2018,braunDynamicsFunctionTear2015}. Periodic boundary conditions are applied.
\begin{align}
\partial_t h  + \partial_x (h\overline{u}) &= -J + P_c(c-1), \label{eq:B1}\\
\overline{u} &=-\frac{h^2}{12}\partial_x p, \label{eq:B2}\\
p &=-\partial_x^2 h, \label{eq:B3}\\
h(\partial_t c +\overline{u}\partial_x c ) &= \text{Pe}_c^{-1}\partial_x (h\partial_x c)+Jc-P_c(c-1)c, \label{eq:B4}\\
h(\partial_t f +\overline{u}\partial_x f ) &= \text{Pe}_f^{-1}\partial_x (h\partial_x f)+Jf-P_c(c-1)f, \label{eq:B5}
\end{align}
The evaporation function $J$ is

\begin{equation} \label{eq:AJ2}
J(x) = \beta \cdot v_b + (a-\beta \cdot v_b) \exp\left[-(x/x_w)^2 / 2 \right].
\end{equation}

The parameters are identical to those in the radial case, and $r$, $r_w$ have simply been replaced by $x$, $x_w$.

\section{Data preprocessing}\label{secA3}
Given a sequence of images $\{X_k\}_{k=1}^n$, we need to stabilize the location of a dark spot over time, since it can move due to flow or movement of the subject's eye. We aim to have the dark spot close to the middle of the region for each image. We achieve this by taking a moving window approach, working backwards in time in order to start when the spot stands out the most.

Let each image $X_k$ be a two-dimensional array indexed by $(i,j) \in \mathbb{Z}^2$, where $i$ denotes the row index and $j$ denotes the column index. For each frame $k$, we define a rectangular subarray $R_k$ centered at pixel coordinates
$\boldsymbol{c_k}=(c_k^{(i)},c_k^{(j)})$. The horizontal and vertical radii are $r_i$ and $r_j$:
\begin{align}
R_k = X_k\! \left\llbracket c_k^{(i)} - r_i : c_k^{(i)} + r_i,\ c_k^{(j)} - r_j : c_k^{(j)} + r_j \right\rrbracket,
\end{align}
where the bracket notation denotes array slicing.

We initialize the procedure by identifying the dark spot center
$c_n$ in the final image $X_n$ and constructing the corresponding
region $R_n$.
To align an earlier frame $X_k$ with the already-aligned region
$R_{k+1}$, we search for a displacement
\[
\boldsymbol{\delta_k} = (\delta_i,\delta_j) \in \{-s,\ldots,s\} \times \{-s,\ldots,s\},
\]
where $s$ is a prescribed search radius, that minimizes the
$\ell^2$-difference between $R_{k+1}$ and the shifted region in frame $k$.
The optimal displacement is defined by
\[
\boldsymbol{\delta_k^{*}}
=
\arg\min_{\delta_i,\delta_j}
\left\|
R_{k+1}
-
X_k\!\left\llbracket
c_{k+1}^{(i)} + \delta_i \pm r_i,\;
c_{k+1}^{(j)} + \delta_j \pm r_j
\right\rrbracket
\right\|_2 .
\]

For stability, we accept the displacement
$\boldsymbol{\delta_k^{*}}$ only if it yields a sufficient reduction in mismatch
relative to the unshifted case. Let
\[
D_{\mathrm{stay}}
=
\left\|
R_{k+1}
-
X_k\!\left\llbracket
c_{k+1}^{(i)} \pm r_i,\;
c_{k+1}^{(j)} \pm r_j
\right\rrbracket
\right\|_2
\]
denote the discrepancy when no displacement is applied. The displacement
$\boldsymbol{\delta_k^{*}}$ is accepted only if
\[
\left\|
R_{k+1} - R_k(\boldsymbol{\delta_k^{*}})
\right\|_2
<
0.95\, D_{\mathrm{stay}} .
\]
Otherwise, we set $\boldsymbol{\delta_k^{*}} = (0,0)$.

Finally, the window center for frame $k$ is updated according to
\[
\boldsymbol{c_k} = \boldsymbol{c_{k+1}} + \boldsymbol{\delta_k^{*}},
\]
and the process is repeated backward in time until all frames have been
aligned.

Before fitting the experimental data to an evaporation function, we apply a smoothing gaussian filter with a standard deviation of $2$ pixels. We also smooth near the rectangle boundary in order to make the data periodic:
\begin{align}
    I_2(x, y) = \frac{1}{4} \left[\tanh\big(k(x - a)\big) - \tanh\big(k(x - b)\big)\right] \cdot \left[\tanh\big(k(y - a)\big) - \tanh\big(k(y - b)\big)\right],
    \end{align}
where $a$ and $b$ define the location and size of the window, and $k$ defines the sharpness of the transition at the boundaries. For our work, we picked $a = -2.6,$ $b = 2.6,$ and $k = 5$.

For each filtered frame $F_k(x,y)$, we construct the image as
\begin{align}
    F^{\text{new}}_k(x, y) = \mu_k + I_2(x, y) \cdot \left({F}_k(x, y) - \mu_k\right).
\end{align}
where $\mu_k$ is the mean value of the filtered image $F_k(x,y)$.
We normalize the filtered data to its largest initial FL intensity value and then interpolate the selected region into $40\times40$ grid for comparison to the model solutions.

\end{appendices}

\section*{Declarations}

\vspace{1em} \noindent \textbf{Conflict of Interest}. The authors declare no competing interests.

\vspace{1em} \noindent \textbf{Funding}. No funding was received to assist with the preparation of this manuscript.

\section*{Data Availability Statement}
The experimental data were collected by our collaborators in a study conducted at Indiana University~\cite{awisi-gyauChangesCornealDetection2019a}. Although the data are not publicly available, interested parties may contact us, and we will seek permission from our collaborators for access.
\bibliography{chenq_cite}

\end{document}